# A new automated strategy for optimizing inclined interplanetary low-thrust trajectories

B. M. Burhani[a], E. Fantino[a], R. Flores[a,b], M. Sanjurjo-Rivo[c]

[a] Department of Aerospace Engineering, Khalifa University of Science and Technology, P.O. Box 127788, Abu Dhabi (United Arab Emirates)
[b] Centre Internacional de Métodes Numèrics en Enginyeria (CIMNE), Gran Capità s/n, 08034 Barcelona (Spain)
[c] Bioengineering and Aerospace Engineering Department, Universidad Carlos III de Madrid, 28911 Leganés, Madrid (Spain)

**Abstract.** This study proposes a new automated strategy for designing and optimizing three-dimensional interplanetary low-thrust (LT) trajectories. The method formulates the design as a hybrid optimal control problem and solves it using a two-step approach. In Step 1, a three-dimensional model based on generalized logarithmic spirals is used with heuristics in combination with a gradient-based solver to perform an automated multi-objective global search of trajectories and optimize for parameters defining the spirals, the launch date, as well as the number, sequence and configuration of the planetary flybys. In Step 2, candidate solutions from Step 1 are refined by further optimization with a direct method. Results show that, compared to similar algorithms based on two-dimensional models, the strategy implemented in Step 1 leads to better estimates of the optimal trajectories, especially when the orbits of the involved bodies are inclined with respect to the ecliptic plane. The proposed approximate method (Step 1) yields better agreement with high-fidelity solutions (Step 2) in terms of launch, flyby and arrival dates, in-plane and out-of-plane average LT accelerations and propellant consumption, leading to improved convergence when the Step 1 trajectories are employed to initiate the search in Step 2.

## 1    Introduction

Nowadays, electric propulsion (EP) is a well-established space technology. It provides significant leverage for decreasing launch mass, increasing payload and extending the spacecraft operational life because it uses smaller amounts of propellant to generate the same total impulse as state-of-the-art chemical thrusters. Additionally, EP can function for extended periods of time, enabling a gradual trajectory modification. As a result, the use of EP ensures wider launch windows and decreases the sensitivity to launch injection errors. However, due to the low propellant mass flow rates and the dependence on the availability of electrical power, EP currently offers modest thrust levels compared to chemical engines.

The performance of EP for deep space scenarios has been proven in important missions, such as NASA Deep Space 1 (DS1) (Rayman et al. 2000), ESA SMART-1 (Racca et al. 2002) and the Dawn rendezvous mission to the asteroid Vesta and the dwarf planet Ceres (Rayman et al. 2006). DS1 (1998), the first interplanetary spacecraft to utilize EP technology, flew by the asteroid Braille and successfully demonstrated the capability to achieve large velocity variations with significant propellant savings. In 2003, ESA proved the viability of EP technology by sending the SMART-1 spacecraft to an orbit around the Moon, consuming just 82 kg of propellant (22% of the launch mass). In 2007, the EP-powered Dawn probe successfully accomplished the highly-challenging objective of entering orbit around Vesta and Ceres, two of the largest objects in the main asteroid belt.

When combined with gravity assists (GAs, or flybys), EP can be very beneficial in terms of time-of-flight reduction and propellant savings, opening up new avenues for solar system exploration (Williams and Coverstone-Carroll 1997). Transfer time and propellant mass fraction are indeed the two most-frequently used measures of the efficiency of an interplanetary trajectory. For this reason, great efforts have been devoted to the development of design methodologies to minimize these two parameters. The conventional approach consists in formulating the trajectory design as an optimization problem, with an objective function to be minimized and a set of constraints to be satisfied. In particular, the optimization process requires the definition of appropriate objectives, the modeling of the system's dynamics and constraints, the development of a solution approach and the selection of a solution technique (Betts 1998; Conway 2012). When planetary gravity assists are incorporated in the trajectory, the complexity of the problem increases because of the additional parameters involved (e.g., the number of flybys, the chronological order of their execution and their geometrical characteristics). A higher degree of complexity also translates into more time-consuming computations. A planetary gravity assist is typically treated as a discrete event. The coupling with the spacecraft's continuous low-thrust (LT) dynamics yields a hybrid optimal control problem that may require mixed-integer programming formulation. In this case, adopting conventional optimal control techniques (such as direct methods) is inefficient due to their poor handling of integer constraints (Morante et al. 2021). In addition, if a high-fidelity dynamical model is adopted for preliminary trajectory design, the solution-finding task becomes impractical (Olympio 2008a; Whiffen 2006). Furthermore, even if the problem can be solved using robust and accurate optimization techniques, a sufficiently good initial guess is usually needed for the search to converge (Conway 2010). Therefore, in the initial stages of mission design, approximate solutions using simplified models are usually preferred. To reduce the dynamical complexity of the problem, the trajectory is often split into a series of patched two-body segments. By assuming that the spheres of influence of the planets are infinitesimal (the so-called zero-radius sphere of influence approximation, ZRSI), the flybys can be modelled as instantaneous changes in the spacecraft's heliocentric velocity (see, e.g., Olympio 2008b). As a result, the trajectory is represented by a series of discrete events (i.e., the flybys) over a continuous path (i.e., the heliocentric transfer). Such combination of discrete events and continuous dynamics transforms the LT multi-gravity-assist trajectory design problem into a hybrid optimal control problem (HOCP). For the description of the general frameworks of HOCPs and their mathematical formulations, the reader may refer to Buss et al. (2002) and Branicky et al. (1998). For HOCP methods specifically tailored for space missions, the works by Ross and D'Souza (2005) and by Chilan and Conway (2013) are recommended.

For a long time (more than six decades to date), a simplified class of hybrid problems known as multi-phase optimal control problems (MOCPs) has been vigorously pursued. In this formulation, the sequence of planetary encounters is *a priori* specified, and MOCPs are regarded as optimal control problems (OCPs) with interior point constraints. For example, the MOCP formulation is implemented in the Solar Electric Propulsion Trajectory Optimization Program by Sauer (1973), which relies on a nonlinear programming (NLP) solution approach to optimize LT multi-gravity assist trajectories. The Gravity Assisted Low-Thrust Local Optimization Program (GALLOP) (McConaghy et al. 2003) is another LT optimization tool based on the MOCP formulation. It combines the Sims-Flanagan direct transcription method (Sims and Flanagan 2000) with a gradient-based NLP solver. In another study, Yam et al. (2011) used the monotonic basin hopping (MBH) technique in conjunction with Sims-Flanagan's transcription to solve continuous-thrust problems with the MOCP formulation.

Because of the inherent combination of integer (discrete) and real (continuous) decision variables, in addition to the nonlinearity and non-convexity of the continuous dynamics (which gives rise to many local optima), determining the optimal sequence of flybys and the optimal LT steering law is challenging. To find globally-optimal trajectories, one simple approach is to form all the possible combinations of targeted flyby bodies, solve each combination as a different MOCP and identify the best solution among the several outcomes. Even if the mission planner's experience can help reduce the size of the search space, this strategy is time consuming. Thus, an algorithm (such as a heuristic method) that searches intelligently among all possibilities, is preferred. Some researchers have been able to develop medium-fidelity tools (which use approximate dynamical models with a good compromise between accuracy and computational speed) that can automatically search for optimal flyby sequences. However, most of these optimizers have been designed for trajectories that involve impulsive maneuvers (e.g., Gad and Abdelkhalik 2011 and Englander et al. 2012). Englander and Conway (2017) developed a tool for LT trajectory design which solves single-objective problems by combining two nested optimization techniques: the outer loop selects the flyby sequence using a genetic algorithm, while the inner loop solves the corresponding MOCP using the MBH global search procedure with Sims-Flanagan's transcription method. The combination yields a medium-fidelity algorithm. Englander et al. (2015) and Vavrina et al. (2015) later extended this methodology to include multi-objective problems and spacecraft system optimization, respectively.

For conceptual trajectory assessment, Sims-Flanagan transcription-based methods are frequently used in most tools seeking to optimize LT multi-flyby trajectories. Depending on the complexity of the problem, computing times for such programs range from several hours to days (Englander and Conway 2017; Englander et al. 2015). Lower-fidelity but faster evaluations can fulfill the purpose of giving quick, comprehensive overviews of the solution space while acting as an adequate starting point for more advanced optimizers. The conventional approach of this kind employs low-fidelity approximations of the LT trajectory, such as shape-based representations, in conjunction with a global search algorithm. Numerous authors have used shape-based techniques to effectively produce solutions very close to the optimal ones. Shape-based models assume an analytical expression for the trajectory, from which the associated thrust profile is derived. Examples include the work by Wall and Conway (2009) who modeled the trajectory as an inverse polynomial, the logarithmic spirals by Bacon (1959) and Tsu (2012), the exponential sinusoid due to Petropoulos and Longuski (2004), the three-dimensional shape-based models incorporating spherical coordinates and pseudo-equinoctial elements by Novak and Vasile (2011) and De Pascale and Vasile (2006), respectively, and the finite Fourier series parametrization due to Taheri and Abdelkhalik (2012). On the basis of the thrust profile of the logarithmic spirals, Roa et al. (2016) discovered an entirely new family of generalized logarithmic spirals, the adaptability of which was later enhanced by including a thrust control parameter in the model (Roa and Peláez 2016a) and by extending it to three dimensions (Roa and Peláez 2016b). Recently, Morante et al. (2019) developed the MOLTO-IT tool, a two-step method that combines planar generalized logarithmic spirals (Roa and Peláez 2016a) with a non-dominated sorting genetic algorithm (NSGA-II) to quickly search for LT interplanetary trajectories confined within the plane containing the Sun and the Earth (the ecliptic plane). Compared with similar preliminary tools such as STOUR-LTGA by Petropoulos and Longuski (2004), MOLTO-IT is faster and generates similar or better guesses, very close to the optimal solutions obtained with more advanced, high-fidelity optimizers.

However, MOLTO-IT suffers from slow convergence rates or even fails to converge to a feasible trajectory when the involved bodies are in moderate or highly-inclined orbits relative to the ecliptic plane. In addition, in these cases, the global search provided by the surrogate model in the first step is not representative of the actual solutions. In such situations, modeling the out-of-plane component of the motion would enhance the accuracy and applicability of the method. This paper aims to extend the MOLTO-IT tool to inclined trajectories by implementing the concept of three-dimensional logarithmic spirals proposed by Roa and Peláez (2016b). In Step 1, a three-dimensional model based on generalized logarithmic spirals is used to perform an automated multi-objective global search of trajectories. In Step 2, as in MOLTO-IT, candidate solutions from Step 1 are refined with a direct method.

This paper is organized as follows. Section 2 summarises the dynamical models for the interplanetary segments and the discrete flyby events. Section 3 formulates the LT interplanetary trajectory optimization problem with variable flyby sequences as a HOCP. The two-step solution strategy is discussed in Sect. 4, while Sect. 5 illustrates the validation of the method through an application to a rendezvous mission to Ceres. Section 6 draws the conclusions.

## 2 Modeling

This section describes the dynamical models for both the continuous and discrete dynamics. The continuous state ($\mathbf{x}$) of a spacecraft is defined by its position ($\mathbf{r}$), velocity ($\mathbf{v}$) and mass (m), i.e., $\mathbf{x} = [\mathbf{r}, \mathbf{v}, m]$. Additionally, the discrete state $q \in \{0,1,2\}$ determines the working condition of the electric engine and the presence/absence of a flyby event. Here, '0' and '1' stand for coasting and thrusting phases of the spacecraft during the interplanetary segments, respectively, while '2' refers to the occurrence of gravity assist maneuvers. Furthermore, the control vector $\mathbf{u}$, which represents the thrust of the LT engine, influences the continuous dynamics.

### 2.1 General dynamical models

Assuming that the motion of the spacecraft is influenced only by the gravitational attraction of the Sun and the thrust acceleration $\mathbf{a}_\mathrm{p}$ of the electric engine, the continuous dynamics can be described by the set of differential equations:

$$\begin{cases} \dot{\mathbf{r}} = \mathbf{v} \\ \dot{\mathbf{v}} = -\frac{\mu}{r^3}\mathbf{r} + \mathbf{a}_\mathrm{p}(\mathbf{x}, q, \mathbf{u}, t) \\ \dot{m} = \dot{m}(\mathbf{x}, q, \mathbf{u}, t), \end{cases} \qquad (1)$$

where $\mu$ is the gravitational parameter of the Sun and $r = \|\mathbf{r}\|$ is the heliocentric distance of the spacecraft. The position ($\mathbf{r}$) and velocity ($\mathbf{v}$) are represented in a global Sun-centered ecliptic reference frame. The thrust acceleration and the mass evolution depend on the on/off state of the thruster, such that

$$\begin{cases} \mathbf{a}_\mathrm{p} = \frac{T}{m}\hat{\mathbf{u}}, \quad \dot{m} = -\frac{T}{I_{sp}g_0} & \text{if} \quad q = 1 \\ \mathbf{a}_\mathrm{p} = \mathbf{0}, \quad\quad \dot{m} = 0 & \text{if} \quad q = 0, \end{cases} \qquad (2)$$

with $T$ being the thrust magnitude and $\hat{\mathbf{u}}$ the unit vector in the direction of thrust (hereafter, hats shall denote unit vectors), while $I_{sp}$ is the specific impulse of the engine and $g_0 = 9.81 \ m/s^2$. $\hat{\mathbf{u}}$ is characterized by the angles $\rho$ and $\beta$ measured in a local reference frame $\{\hat{\mathbf{r}}, \hat{\mathbf{w}}, \hat{\mathbf{s}}\}$ with origin at the spacecraft (see Fig. 1):

$$\hat{\mathbf{r}} = \frac{\mathbf{r}}{\|\mathbf{r}\|}, \qquad \hat{\mathbf{w}} = \frac{\mathbf{r} \times \mathbf{v}}{\|\mathbf{r} \times \mathbf{v}\|}, \qquad \hat{\mathbf{s}} = \hat{\mathbf{w}} \times \hat{\mathbf{r}}. \qquad (3)$$

$\rho$ is the angle from $\hat{\mathbf{s}}$ to the projection of $\hat{\mathbf{u}}$ on the orbital plane, measured positively in the counterclockwise direction with respect to $\hat{\mathbf{w}}$, while $\beta$ denotes the angle between $\hat{\mathbf{w}}$ and $\hat{\mathbf{u}}$. Hence, the thrust direction is computed as,

$$\hat{\mathbf{u}} = \sin\beta \sin\rho \; \hat{\mathbf{r}} + \sin\beta \cos\rho \; \hat{\mathbf{s}} + \cos\beta \; \hat{\mathbf{w}}. \qquad (4)$$

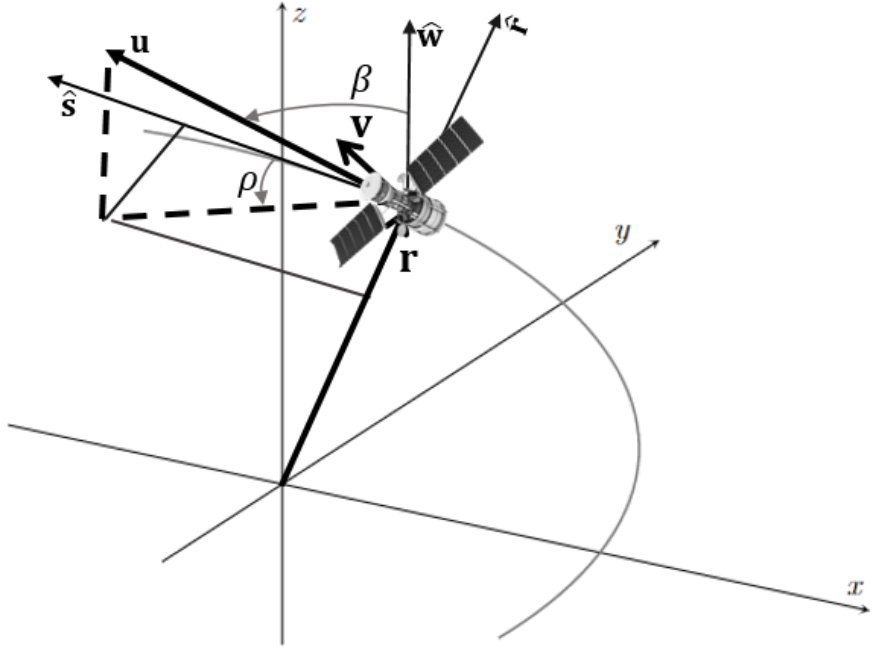

**Fig. 1.** Thrust vector and the associated thrust angles for the continuous dynamics.

In the case of the discrete dynamics associated with flyby events, the ZRSI approximation is adopted to model the planetary gravity assists. It is assumed that the heliocentric position and the mass of the spacecraft do not change throughout the flyby maneuver, which is considered unpowered. The heliocentric post-flyby velocity $\mathbf{v}^+$ can be calculated as a function of the pre-flyby velocity $\mathbf{v}^-$, the planet's heliocentric velocity $\mathbf{v}_b$, the pericenter radius $r_p$, and the B-plane angle $\zeta$. The B-plane is orthogonal to $\mathbf{v}_\infty{}^-$ (the planetocentric inbound velocity) and contains the center of the planet (Greenberg et al. 1988; Farnocchia et al. 2019), see Fig. 2. $\zeta$ is the angle (measured clockwise inside the B-plane when looking in the direction of $\mathbf{v}_\infty{}^-$) between $\hat{\mathbf{f}}$ and the flyby plane, where $\hat{\mathbf{f}}$ is contained on the B-plane and is parallel to the ecliptic. The orientation of $\hat{\mathbf{f}}$ is chosen such that $\hat{\mathbf{f}} \cdot \mathbf{v}_b < 0$.

Assuming that the gravitational parameter of the flyby planet is given by $\mu_b$, the maneuver is described as:

$$\begin{cases} \mathbf{v}_\infty{}^- = \mathbf{v}^- - \mathbf{v}_b \\ \delta = 2\sin^{-1}\left(\frac{1}{1+r_p v_\infty{}^2/\mu_b}\right) \\ \mathbf{v}_\infty{}^+ = v_\infty\left(\cos\delta\,\hat{\mathbf{i}} + \cos\zeta\sin\delta\,\hat{\mathbf{j}} + \sin\zeta\sin\delta\,\hat{\mathbf{k}}\right) \\ \mathbf{v}^+ = \mathbf{v}_\infty{}^+ + \mathbf{v}_b \end{cases} \quad \text{if } q = 2. \tag{5}$$

In Eq. 5, $v_\infty = \|\mathbf{v}_\infty{}^-\| = \|\mathbf{v}_\infty{}^+\|$ is the flyby hyperbolic excess speed, and the unit vectors are defined as:

$$\hat{\mathbf{i}} = \frac{\mathbf{v}_\infty{}^-}{\|\mathbf{v}_\infty{}^-\|}, \qquad \hat{\mathbf{j}} = \frac{\mathbf{i}\times\mathbf{v}_b}{\|\mathbf{i}\times\mathbf{v}_b\|}, \qquad \hat{\mathbf{k}} = \hat{\mathbf{i}} \times \hat{\mathbf{j}}. \tag{6}$$

Here, $\mathbf{v}_\infty{}^+$ stands for the outbound spacecraft velocity relative to the planet, whereas $\delta$ is the deflection angle on the flyby plane.

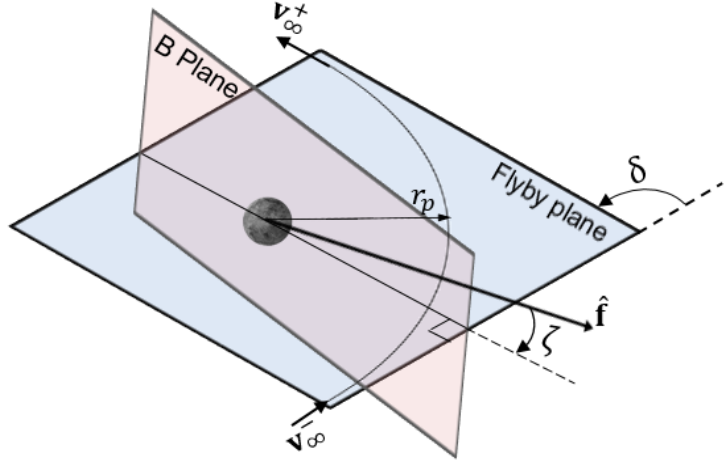

**Fig. 2.** The flyby geometry and the B-plane.

### 2.2 Simplified dynamical model: 3D generalized logarithmic spirals

In the preliminary stages of trajectory design, one can neglect complex constraints and adopt a simplified dynamical model to reduce the computational burden. In our case, shape-based methods constrain the profile of the steering law.

The simplified dynamical model employed in this study was proposed by Roa and Peláez (2016b). During the propelled phases, the motion of the spacecraft is described by three-dimensional generalized logarithmic spirals. In this formulation, the trajectory is split into in-plane and out-of-plane components. The former defines the base spiral (gray in Fig. 3) and includes a set of spiral control parameters, which can be adjusted to fulfill the trajectory constraints. There are multiple options for the out-of-plane motion, such as helixes, polynomial shaping laws, and Fourier series (Roa and Peláez 2016b). The main advantage of this formulation is that it allows expressing the trajectory in analytical form in terms of a small set of parameters, substantially reducing the complexity of the problem.

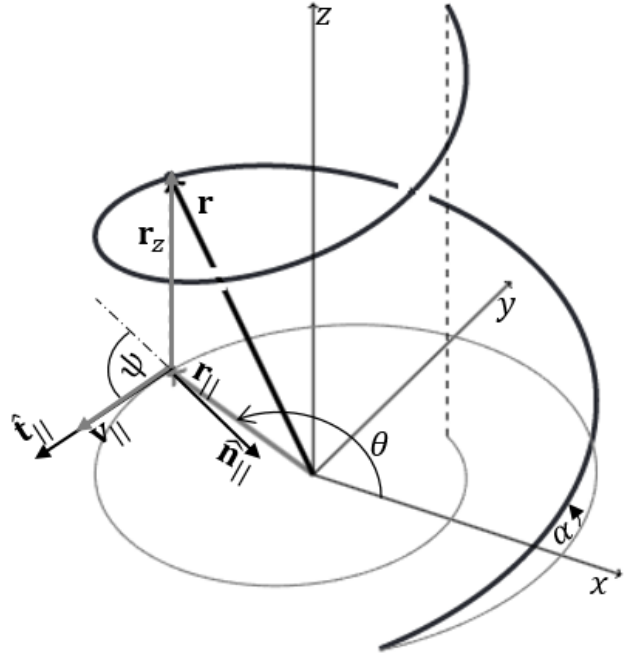

**Fig. 3.** Geometry of a 3D generalized logarithmic spiral.

The acceleration due to the thrust $\mathbf{a}_p$ is split into ecliptic $\mathbf{a}_{p,\parallel}$ (which we will refer to as in-plane) and normal-to-the-ecliptic $\mathbf{a}_{p,z}$ (out-of-plane) components,

$$\mathbf{a}_p = \mathbf{a}_{p,\parallel} + \mathbf{a}_{p,z} \ . \tag{7}$$

Similarly, position and velocity vectors are expressed as

$$\mathbf{r} = \ \mathbf{r}_{\parallel} + \mathbf{r}_z, \tag{8}$$

$$\mathbf{v} = \ \mathbf{v}_{\parallel} + \mathbf{v}_z. \tag{9}$$

**Planar motion.** Following Roa and Peláez (2016a and 2016b), we use logarithmic spirals for the planar component of thrust

$$\mathbf{a}_{p,\parallel} = \mu \left[ \left( \frac{\xi - 1}{r_{\parallel}^{3}} + \frac{1}{r^3} \right) r_{\parallel} \cos\psi \ \hat{\mathbf{t}}_{\parallel} + \left( \frac{2(1-\xi)}{r_{\parallel}^{3}} - \frac{1}{r^3} \right) r_{\parallel} \sin\psi \ \hat{\mathbf{n}}_{\parallel} \right], \tag{10}$$

where $\xi$ is a spiral control parameter and $\psi$ is the planar flight path angle, defined by (Fig. 3)

$$\cos\psi = \frac{\mathbf{r}_{||}\cdot\mathbf{v}_{||}}{r_{||}v_{||}},$$  (11)

with $r_{||} = \|\mathbf{r}_{||}\|$ and $v_{||} = \|\mathbf{v}_{||}\|$. The in-plane tangential $\hat{\mathbf{t}}_{||}$ vector is given by

$$\hat{\mathbf{t}}_{||} = \frac{\mathbf{v}_{||}}{\|\mathbf{v}_{||}\|},$$  (12)

and the in-plane normal $\hat{\mathbf{n}}_{||}$ is perpendicular to $\hat{\mathbf{t}}_{||}$ and points towards the center of curvature. The base spiral, which describes the planar motion, is obtained substituting the planar thrust acceleration (Eq. 10) into Eq. 1, ignoring the evolution of mass and out-of-plane states. Thus, the continuous dynamical model becomes (Roa and Peláez 2016a):

$$\begin{cases} \frac{dv_{||}}{dt} = \mu \frac{\xi-1}{r_{||}^2}\cos\psi \\ \frac{dr_{||}}{dt} = v_{||}\cos\psi \\ \frac{d\theta}{dt} = \frac{v_{||}}{r_{||}}\sin\psi \\ \frac{d\psi}{dt} = \frac{\sin\psi}{v_{||}r_{||}^2}\left(2(1-\xi) - r_{||}v_{||}^2\right) \end{cases}.$$  (13)

Eq. 13 can be integrated analytically, yielding two constants of motion, $K_1$ and $K_2$, determined by the initial conditions (henceforth, the subscript 0 denotes initial values):

$$K_1 = v_{0,||}^2 - \frac{2\mu}{r_{0,||}}(1-\xi),$$  (14)

$$K_2 = r_{0,||}v_{0,||}^2\sin\psi_0.$$  (15)

$K_1$ and $K_2$ are reminiscent of specific energy and specific angular momentum, respectively Roa and Peláez (2016b).

In a similar way as for Keplerian orbits, the sign of $K_1$ discriminates three families of spirals: elliptic ($K_1 < 0$), parabolic ($K_1 = 0$) and hyperbolic ($K_1 > 0$). Hyperbolic spirals are further classified as either type I or type II depending on whether $K_2$ is smaller or larger than $2(1-\xi)$. The analytical expressions for the base spirals in terms of the initial conditions, the constants of motion, the spiral anomaly $\beta(\theta)$ and the spiral control parameter are summarized in Table 1. The angle $\beta$ is a measure of the angular displacement from the periapsis (in the case of elliptic and hyperbolic type II spirals) or from the asymptote (in the case of hyperbolic type I spirals) of the base spiral (see Roa and Peláez (2016b) for more details).

**Table 1.** Analytical expressions for the base spirals.

| Family of base spiral | In-plane component of the trajectory |
|---|---|
| Elliptic | $r_{\parallel} = \dfrac{-2(1-\xi) + K_2}{(K_1 + K_1 K_2 \cosh \beta(\theta))}$ |
| Parabolic | $r_{\parallel} = r_{0,\parallel} e^{(\theta - \theta_0) \cot \psi}$ |
| Hyperbolic Type I | $r_{\parallel} = \dfrac{k_1}{\sinh(\beta(\theta)/2) \left[ k_2 \sinh\left(\dfrac{\beta(\theta)}{2}\right) + k_3 \cosh\left(\dfrac{\beta(\theta)}{2}\right) \right]}$ |
| Hyperbolic Type II | $r_{\parallel} = \dfrac{K_2{}^2 - 4(1-\xi)^2}{2K_1(1-\xi) + K_1 K_2 \cos \beta(\theta)}$ |

In Table 1, the expression of $r_{\parallel}$ for the hyperbolic type I spiral contains the quantities $k_1$, $k_2$, and $k_3$, which are given by

$$\begin{cases} k_1 = \varrho l^2 / K_1 \\ k_2 = 4 \varrho (1 - \xi) \\ k_3 = \varrho^2 - K_2^2 \end{cases} \tag{16}$$

with

$$\begin{cases} \varrho = 2(1 - \xi) + l \\ l = \sqrt{4(1 - \xi^2) - K_2^2} \end{cases} \tag{17}$$

The in-plane speed can be obtained from the definition of $K_1$ in Eq. 14 as,

$$v_{\parallel} = \sqrt{K_1 + \frac{2\mu}{r_{\parallel}} (1 - \xi)} . \tag{18}$$

**Out-of-plane motion.** In this work, the polynomial shaping law describes the out-of-plane component of motion, as it leads to better estimation of trajectories compared to helixes and Fourier series. The out-of-plane component $r_z$ of the position is given as a function of the polar angle of the base spiral $\theta$ as

$$r_z = \sum_{i=0}^{n} c_i \theta^i \,, \tag{19}$$

where $c_i$ are the coefficients of the polynomial function of order $n$. Hence, the out-of-plane velocity component can be written as,

$$v_z = \frac{\mathrm{d}r_z}{\mathrm{d}t} = \frac{\mathrm{d}r_z}{\mathrm{d}\theta}\frac{\mathrm{d}\theta}{\mathrm{d}t} = \frac{v_{||}\sin\psi}{r_{||}}\sum_{i=1}^{n} i c_i \theta^{i-1} \,. \tag{20}$$

Equation 16 yields the polynomial steering law for the climb angle $\alpha$ (see Fig. 3):

$$\tan\alpha = \frac{v_z}{v_{||}} = \frac{\sin\psi}{r_{||}}\sum_{i=1}^{n} i c_i \theta^{i-1} \,. \tag{21}$$

Similarly, the out-of-plane acceleration can be computed as:

$$a_z = \frac{\mathrm{d}v_z}{\mathrm{d}t} = \frac{\mathrm{d}}{\mathrm{d}t}\left(v_{||}\tan\alpha\right) = v_{||}\frac{\mathrm{d}}{\mathrm{d}t}(\tan\alpha) + \tan\alpha\frac{\mathrm{d}v_{||}}{\mathrm{d}t} =$$
$$v_{||}\frac{\mathrm{d}}{\mathrm{d}t}(\tan\alpha) + \tan\alpha\left(\mu\frac{\xi-1}{r_{||}^2}\cos\psi\right) \tag{22}$$

To determine the derivative of $\tan\alpha$ in Eq. 21, we differentiate Eq. 20 with respect to time, yielding:

$$a_z = \frac{\mathrm{d}v_z}{\mathrm{d}t} = \frac{\mathrm{d}v_z}{\mathrm{d}\theta}\frac{\mathrm{d}\theta}{\mathrm{d}t} = \frac{v_{||}\sin\psi}{r_{||}}\frac{\mathrm{d}}{\mathrm{d}\theta}\left(\frac{v_{||}\sin\psi}{r_{||}}\sum_{i=1}^{n} i c_i \theta^{i-1}\right). \tag{23}$$

Therefore,

$$a_z = \frac{v_{||}\sin\psi}{r_{||}}\left(\frac{v_{||}\sin\psi}{r_{||}}\sum_{i=2}^{n} i(i-1)c_i\theta^{i-2} + \frac{\mathrm{d}}{\mathrm{d}\theta}\left(\frac{v_{||}\sin\psi}{r_{||}}\right)\sum_{i=1}^{n} i c_i\theta^{i-1}\right). \tag{24}$$

Simplifying Eq. 24 and making use of Eqs. 13 and 14 gives

$$a_z = \frac{v_{||}^2\sin^2\psi}{r_{||}^2}\sum_{i=2}^{n} i(i-1)c_i\theta^{i-2} - \frac{\sin\psi\cos\psi}{r_{||}^2}\left(K_1 + v_{||}^2\right)\sum_{i=1}^{n} i c_i\theta^{i-1} +$$
$$\mu\frac{\xi-1}{r_{||}^2}\cos\psi\frac{\sin\psi}{r_{||}}\sum_{i=1}^{n} i c_i\theta^{i-1}. \tag{25}$$

Using Eq. 21, the above expression for $a_z$ can be rewritten as:

$$a_z = v_{||}\left(\frac{v_{||}\sin^2\psi}{r_{||}^2}\sum_{i=2}^{n} i(i-1)c_i\theta^{i-2} - \frac{\sin\psi\cos\psi}{v_{||}r_{||}^2}\left(K_1 + v_{||}^2\right)\sum_{i=1}^{n} i c_i\theta^{i-1}\right) +$$
$$\mu\frac{\xi-1}{r_{||}^2}\cos\psi\tan\alpha. \tag{26}$$

Hence, by comparing Eqs. 22 and 26, one can obtain the time derivative of $\tan\alpha$ as,

$$\frac{\mathrm{d}}{\mathrm{d}t}(\tan\alpha) = \frac{v_{||}\sin^2\psi}{r_{||}^2}\sum_{i=2}^n i(i-1)c_i\theta^{i-2} - \frac{\sin\psi\cos\psi}{v_{||}r_{||}^2}(K_1 + v_{||}^2)\sum_{i=1}^n ic_i\theta^{i-1}. \quad (27)$$

The out-of-plane acceleration has gravitational $(a_{g,z})$ and thrust $(a_{p,z})$ components, i.e.,

$$a_z = a_{g,z} + a_{p,z}. \quad (28)$$

Hence,

$$v_{||}\frac{\mathrm{d}}{\mathrm{d}t}(\tan\alpha) + \tan\alpha\left(\mu\frac{\xi-1}{r_{||}^2}\cos\psi\right) = -\mu\frac{r_z}{r^3} + a_{p,z}, \quad (29)$$

yielding

$$a_{p,z} = \mu\left(\frac{\xi-1}{r_{||}^2}\right)\tan\alpha\cos\psi + v_{||}\frac{d}{dt}(\tan\alpha) + \mu\frac{r_z}{r^3}. \quad (30)$$

The magnitude of in-plane and out-of-plane thrust accelerations can be used to determine the total impulse $\Delta V$ along the associated spiral as,

$$\Delta V = \int_{t_1}^{t_2}\sqrt{(a_{p,||}^2 + a_{p,z}^2)}\,dt, \quad (31)$$

where $t_1$ and $t_2$ represent the times at which the spiral motion starts and ends, respectively. The total impulse is then used to calculate the propellant mass fraction using the rocket equation, assuming a constant specific impulse:

$$\frac{m_p}{m_0} = 1 - e^{-\frac{\Delta V}{I_{sp}g_0}}, \quad (32)$$

where $m_p$ is the propellant mass consumed and $m_0$ is the initial mass of the spacecraft.

## 3    Hybrid optimal control problem

In this section, we describe the approach used to write the optimal design of three-dimensional LT interplanetary trajectories as a hybrid optimal control problem (HOCP). Several mathematical frameworks for HOCP formulation exist in the literature (Buss et al. 2002; Branicky et al. 1998; Ross and D'Souza 2005; Chilan and Conway 2013). In this paper, for the sake of simplicity, we adopt the one based on the work of Buss et al. (2002). Let us assume the discrete flyby events at transition times $t_j$ are given by $n_s$ discontinuity surfaces such that $s_j = 0$ for $j = \{1, 2, \ldots, n_s\}$ (see Fig. 4).

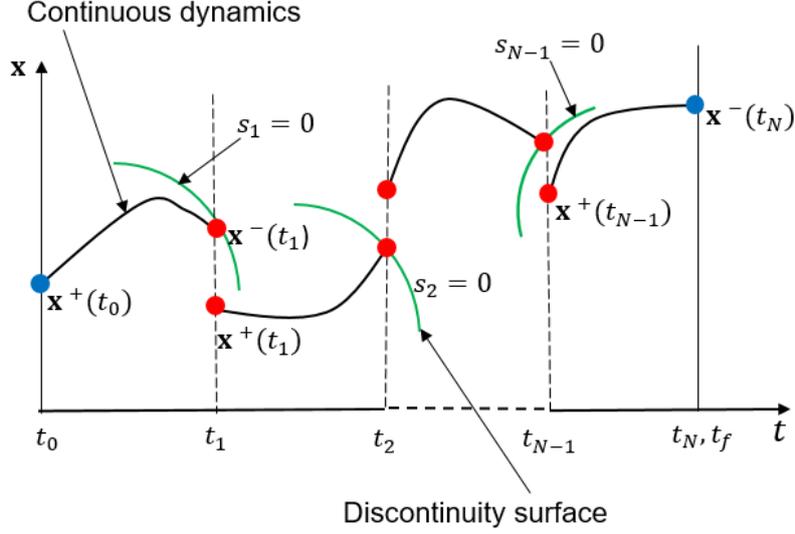

**Fig. 4.** Illustration of a hybrid dynamical system.

Thus, the evolution of the state vector with time is given by:

$$\dot{\mathbf{x}}(t) = f(\mathbf{x}, q, \mathbf{u}, t) \quad \text{if} \quad s_j(\mathbf{x}, q, \mathbf{u}, t) \neq 0 \; \forall j, \tag{33}$$

and

$$[\mathbf{x}^+(t_j), q^+(t_j)] = \varphi_j(\mathbf{x}^-(t_j), q^-(t_j), \mathbf{u}, t_j) \quad \text{otherwise}, \tag{34}$$

where $f$ is a function which governs the continuous dynamics (Eq. 1), while $\varphi$ is a set of transition map functions that describes the discontinuous behavior characterized by flybys (Eq. 5). In Eq. 34, $\mathbf{x}^-$ and $q^-$ stand for continuous and discrete states just before a GA maneuver at $t = t_j$; $\mathbf{x}^+$ and $q^+$ represent the corresponding states immediately after. Suppose that there are $N$ interplanetary legs, and hence $N - 1$ flyby events. Their time-ordered sequence can be expressed as

$$\sigma = [t_1, \ldots, t_{N-1}], \tag{35}$$

where at each of the flyby times the following condition must hold

$$s_i = \|\mathbf{r}(t_i) - \mathbf{r}_{b_i}(t_i)\| = 0 \quad \text{for} \quad i = 1, \ldots, N-1, \tag{36}$$

$\mathbf{r}_{b_i}$ denoting the position of the $i$-th flyby body. Note that number of discontinuities, the sequence of planets and the times of the flybys are unknown, and should be determined as a part of the optimal solution.

Generally, the hybrid optimal control problem involves finding the optimal values of the control ($\mathbf{u}$) and the associated continuous ($\mathbf{x}$) and discrete ($q$) states that minimize a certain objective function, $J$, subject to the constraints of the problem. Mathematically, it can be written as

Minimize:

$$J = \Phi(\mathbf{x}^+(t_0), \dots, \mathbf{x}^-(t_N); \; q^+(t_0), \dots, q^-(t_N); \; t_0, \dots, t_N) \; + \int_{t_0}^{t_f} \Psi(\mathbf{x}, q, \mathbf{u}, t) \, \mathrm{d}t, \quad (37)$$

subject to dynamical constraints (Eqs. 33-34), with

$$\mathbf{x}(t) \; \in X, \tag{38}$$

$$q(t) \; \in Q, \tag{39}$$

$$\mathbf{u}(t) \; \in U, \tag{40}$$

$$g_l \leq g(\mathbf{x}, q, \mathbf{u}, t) \leq g_u, \tag{41}$$

$$h(\mathbf{x}(t_0), q(t_0), t_0, \mathbf{x}(t_f), q(t_f), t_f) \leq 0, \tag{42}$$

$$t \in [t_0, t_f]. \tag{43}$$

The objective function $J$ is comprised of two parts, namely the Mayer and Bolza terms. The Mayer term is a function $\Phi$ of discrete event times $t_0, \dots, t_N$, as well as associated states just before and after the corresponding discrete events. The Bolza term contains an integrand, $\Psi$, which is a function of continuous states, controls, and time. The constraints of the problem include the dynamical constraints (Eqs. 33-34), the admissible set of continuous and discrete states $X$ and $Q$ (Eqs. 38-39), admissible controls $U$ (Eq. 40), path constraints (Eq. 41), boundary conditions (Eq. 42), and time bounds (Eq. 43). In Eq. 41, $g_l$ and $g_u$ represent the minimum and maximum values of the path constraint function, while in Eqs. 37, 42, and 43, $t_0$ and $t_f$ stand for the initial and final dates.

In this study, instead of trying to find a unique optimal solution as in single-objective problems, we perform a multi-objective optimization, obtaining a set of candidate solutions. Hence, the HOCP seeks to determine the best number and sequence of flybys as well the LT steering law that minimizes a set of performance indices $J = [J_1, ..., J_{n_j}]$, where $n_j$ is the number of objective functions. However, the high-fidelity form of the HOCP presented above is far too complex to handle with a reasonable amount of computer resources. As a result, simplifications in the dynamics and/or constraints are commonly assumed.

## 4    Solution approach

In this section, we present a solution approach for solving the HOCP which builds on the work of Morante et al. (2019), who introduced the MOLTO-IT algorithm. The MOLTO-IT strategy involves two steps. The first (Step 1 hereafter) is a broad search of optimal trajectories which determines the number and sequence of flybys. In this step, a simplified dynamical model is used; the continuous dynamics is approximated with planar generalized logarithmic spirals (Eq. 13) and a predefined control law, while the discrete dynamics associated with GA maneuvers is described with a two-dimensional model. The algorithm returns a set of optimal trajectories based on two performance criteria: transfer time and propellant mass consumption. It also provides information on the mission design variables, including launch, flyby, and arrival dates, flyby sequence as well as the launch hyperbolic excess speed. The user can then select the solution that strikes the best compromise between mission requirements. The selected trajectory is used as an initial guess in the second step, named Step 2. In Step 2, a more accurate 3D dynamical model is used (Eqs. 1 and 5). In addition, the thrust profile is determined automatically by a high-fidelity NLP solver, which enforces the complete set of constraints.

In this study, a similar approach is adopted. Rather than using a 2D model in Step 1, this contribution seeks to incorporate the out-of-plane component of motion. This improved algorithm, 3D MOLTO-IT, relies on a simplified 3D dynamical model based on logarithmic spirals connected with unpowered GAs for step 1. The details of both steps are described next.

### 4.1 Step 1: Global multi-objective search

Using a predefined control law, 3D MOLTO-IT Step 1 reduces the infinite-dimensional HOCP to a mixed integer nonlinear programming problem (MINLP) consisting of a small number of continuous (real) and discrete (integer) optimization parameters. In each interplanetary leg, a thrust-coast-thrust sequence is assumed if the target is to rendezvous with the end planet of the leg (see Fig. 5). If the goal is to perform a flyby with the planet, a thrust-coast sequence is adopted instead (see Fig. 6). In the rendezvous case, it is not possible to match the position and velocity of the second planet simultaneously with just one spiral arc. In general, the initial and terminal points of the interplanetary legs have different values of the constants of motion $K_1$ and $K_2$. However, only $K_1$ is affected by the spiral control parameter (see Eqs. 14-15). As a result, an intermediate coast (Keplerian) arc is introduced in order to connect the departure and arrival spirals with appropriate values of $K_2$. Additional coast and thrust segments may be added for flexibility. However, this increases the computational load significantly. When the target is to flyby with the second planet, there is no need for the second spiral as only the final position constraint must be enforced. Consequently, a thrust-coast sequence is assumed, with the coast arc providing flexibility without significant computational burden. Nevertheless, the position of the switching points between the first spiral and the coast arc, $\theta_A$, and between the coast arc and the second spiral, $\theta_B$, are unknown (see Figs. 5 and 6). They can be optimized to satisfy the leg constraints minimizing the mission cost. Note that, in Figs. 5 and 6, $\xi_1$ and $\xi_2$ denote the spiral control parameters of the respective arcs, while $t_0$ and $t_F$ represent the departure and arrival times of the interplanetary leg.

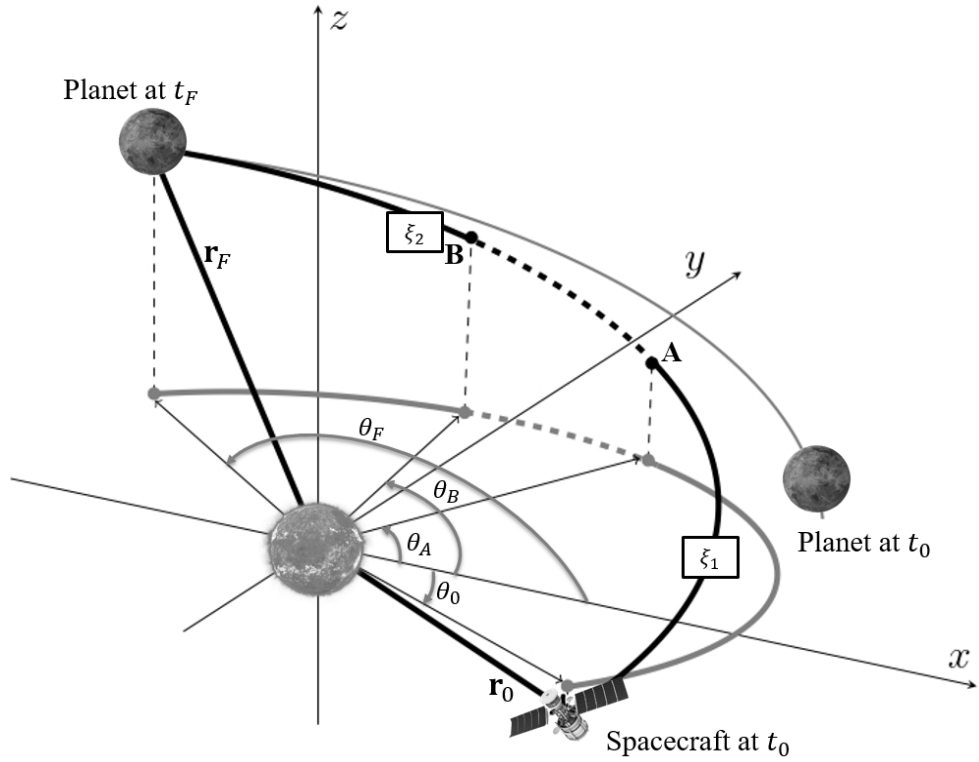

**Fig. 5.** Rendezvous leg configuration.

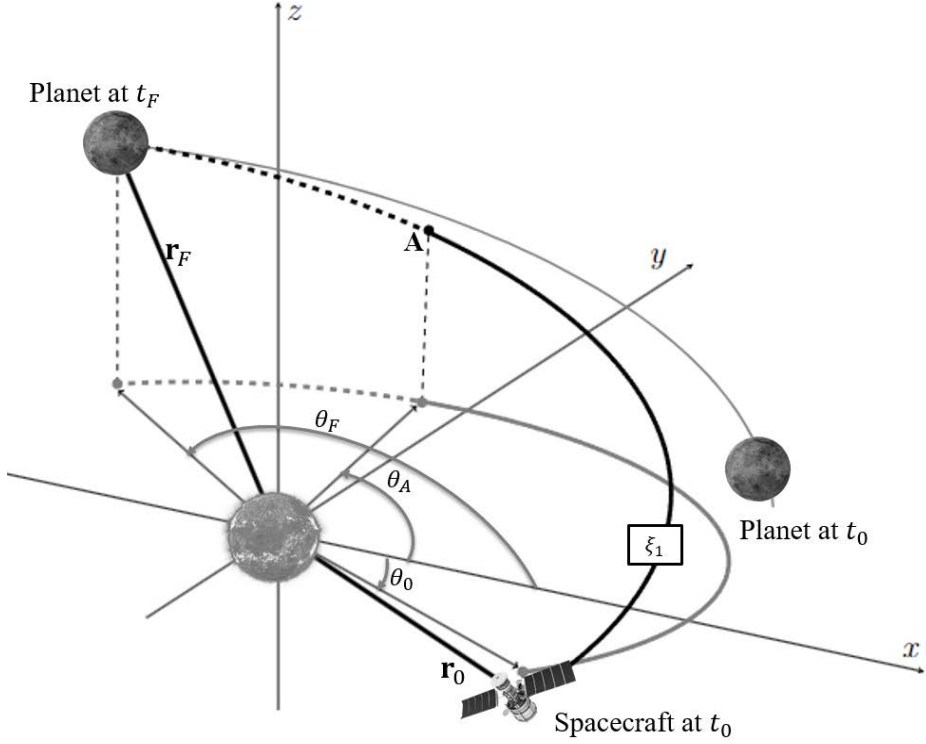

**Fig. 6.** Flyby leg configuration.

Besides the boundary conditions associated with the departure and arrival points of each leg, additional constraints arise from the out-of-plane motion. For the polynomial shaping law, the function must have at least four coefficients Roa and Peláez (2016b). Additional coefficients may be added for flexibility and better trajectory estimation at the cost of computational effort. A fourth-order polynomial has been selected for the 3D MOLTO-IT algorithm, since it offers a good compromise between flexibility and computational burden. The out-of-plane components of position $r_z(\theta)$ and velocity $v_z(\theta)$ are given as (see Eqs. 19-20),

$$r_z(\theta) = c_0 + c_1\theta + c_2\theta^2 + c_3\theta^3 + c_4\theta^4, \tag{44}$$

$$v_z(\theta) = \frac{v_{||}\sin\psi}{r_{||}}(c_1 + 2c_2\theta + 3c_3\theta^2 + 4c_4\theta^3). \tag{45}$$

In this study, the first two coefficients, $c_0$ and $c_1$, are computed from the initial out-of-plane position and velocity and the remaining coefficients $c_2$, $c_3$ and $c_4$.

$$c_1 = \frac{r_{\parallel,0}}{v_{\parallel,0}\sin\psi_0} v_z(\theta_0) - \left(2c_2\theta_0 + 3c_3\theta_0{}^2 + 4c_4\theta_0{}^3\right), \tag{46}$$

$$c_0 = r_z(\theta_0) - \left(c_1\theta_0 + c_2\theta_0{}^2 + c_3\theta_0{}^3 + c_4\theta_0{}^4\right). \tag{47}$$

The treatment of the remaining coefficients depends on the spiral under consideration, i.e., the first or the second. For the second spiral, they are considered optimization variables such that the final out-of-plane states match those of the second planet, and the total $\Delta V$ along the spiral is minimized. In the case of the first spiral, the final out-of-plane position and velocity are unknown, requiring a different set of constraints. We assume that the out-of-plane motion is not significant in the first spiral, meaning that almost all of the out-of-plane impulse is applied in the second spiral. The reason is that the magnitudes of out-of-plane position and velocity are small compared with their in-plane counterparts. Thus, the LT engine typically attempts to correct the in-plane components first. Consequently, the out-of-plane accelerations corresponding to the initial ($\theta = \theta_0$), middle ($\theta = (\theta_0 + \theta_A)/2$), and final ($\theta = \theta_A$) points of the first spiral are set to zero:

$$\begin{cases} a_{p,z}(\theta_0) = 0 \\ a_{p,z}((\theta_0 + \theta_A)/2) = 0 \\ a_{p,z}(\theta_A) = 0. \end{cases} \tag{48}$$

Thus, for the first spiral, the values of coefficients $c_2$, $c_3$ and $c_4$ are obtained enforcing Eq. 48.

Due to the mix of real and integer variables and the requirement to evaluate many different scenarios at the same time, population-based heuristic algorithms, such as genetic algorithms or particle swarm optimization, are suitable solution techniques. They struggle, however, with the nonlinear constraints. Gradient-based solvers can deal with such constraints more effectively. As a result, two nested optimization schemes are used in the 3D tool (the same strategy adopted by MOLTO-IT); a heuristic algorithm (outer loop) and a gradient-based solver (inner loop). Figure 7 illustrates the general scheme of the 3D MOLTO-IT Step 1. The Appendix presents the meta-code for the heuristic (Algorithm 2 hereafter) and gradient-based (Algorithm 1) schemes.

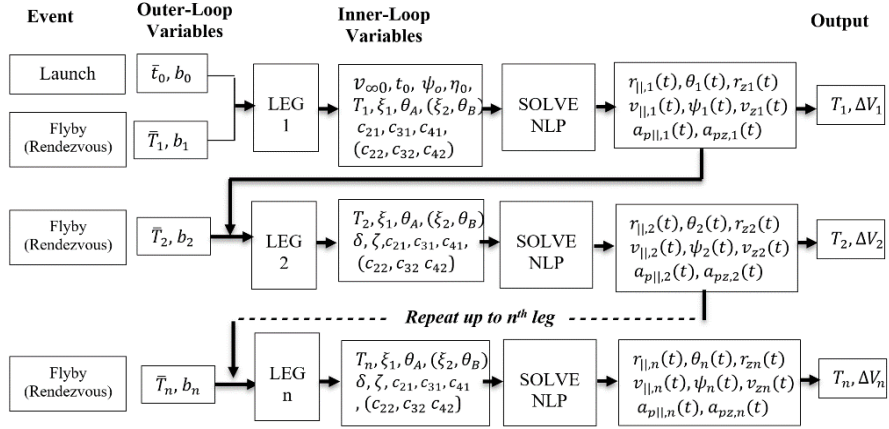

**Fig. 7.** 3D MOLTO-IT Step 1 algorithm scheme. The quantities in parenthesis correspond to a rendezvous event.

**Table 2.** Outer loop (heuristic) algorithm variables

| Variable | Meaning | Lower bound | Upper bound |
|:---:|:---:|:---:|:---:|
| $\bar{t}_0$ | Departure epoch guess | $t_{0\,\min}$ | $t_{0\,\max}$ |
| $\bar{T}_i$ | Leg transfer time | $T_{i\,\min}$ | $T_{i\,\max}$ |
| $b_i$ | Flyby body | - | - |

**Table 3.** Inner loop (NLP) algorithm variables.
Key: L = Launch, F = Flyby, R= Rendezvous

| Variable | Meaning | Lower bound | Upper bound | Initial guess | Event 1 | Event 2 |
|:---:|:---:|:---:|:---:|:---:|:---:|:---:|
| $t_0$ | Departure epoch | $0.9\bar{t}_0$ | $1.1\bar{t}_0$ | $\bar{t}_0$ | L | F or R |
| $v_{\infty 0}$ | Launch excess velocity | $v_{\infty 0,min}$ | $v_{\infty 0,max}$ | $v_{\infty 0,max}$ | L | F or R |
| $\psi_o$ | In-plane launch angle (rad) | $-\pi$ | $-\pi$ | $-\dfrac{\pi}{2}, 0, \dfrac{\pi}{2}$ | L | F or R |
| $\eta_0$ | Off-plane launch angle (rad) | $-\pi/2$ | $\pi/2$ | $0$ | L | F or R |
| $T$ | Leg transfer time | $0.9\bar{T}_i$ | $1.1\bar{T}_i$ | $\bar{T}_i$ | F or R | F or R |

| $\delta$ | Flyby deflection angle (rad) | $-\pi$ | $\pi$ | 0 | F | F or R |
|---|---|---|---|---|---|---|
| $\zeta$ | B-plane angle (rad) | $-\pi/2$ | $\pi/2$ | 0 | F | F or R |
| $\xi_1$ | First spiral control parameter | 0 | 1 | 0.4 | F or R | F or R |
| $\xi_2$ | Second spiral control parameter | 0 | 1 | 0.4 | F or R | R |
| $\theta_A$ | Thrust-to-coast switch angle | 0 | 1 | 0.01 | F or R | F or R |
| $\theta_B$ | Coast-to-thrust switch angle | 0 | 1 | 0.99 | F or R | R |
| $c_{21}$ | First spiral $c_2$ (Eqs. 44-45) | -10 | 10 | 0 | F or R | F or R |
| $c_{31}$ | First spiral $c_3$ | -10 | 10 | 0 | F or R | F or R |
| $c_{41}$ | First spiral $c_4$ | -10 | 10 | 0 | F or R | F or R |
| $c_{22}$ | Second spiral $c_2$ | -10 | 10 | 0 | F or R | R |
| $c_{32}$ | Second spiral $c_3$ | -10 | 10 | 0 | F or R | R |
| $c_{42}$ | Second spiral $c_4$ | -10 | 10 | 0 | F or R | R |

In Algorithm 1, for each interplanetary segment, the algorithm determines first whether the departure point (*Event-1*) represents the launch or a flyby event and also whether the arrival point (*Event-2*) corresponds to a flyby or a rendezvous event. Next, we determine the initial state of the spacecraft from the ephemeris of the departure planet and compute the launch speed of the spacecraft if *Event-1* is the launch scenario. On the other hand, if *Event-1* is characterized by a flyby maneuver, we use the pre-flyby state of the spacecraft and the velocity of the planet as constraints to compute the post-flyby heliocentric velocity of the probe. This information is then used to propagate the trajectory within the first spiral arc to the intermediate point $\theta_A$, where it connects with the subsequent coast arc. Next, the out-of-plane acceleration at both ends and center of the first spiral are set to zero (Eq. 48). The trajectory is then propagated to the end of the coast arc using the state at $\theta_A$. If *Event-2* is a flyby, the position of the spacecraft at the end of the coast arc is equated to that of the planet. If *Event-2* is a rendezvous scenario, the end of the Keplerian segment is used to propagate the trajectory through the second spiral arc to determine the state at the arrival point, $\theta_F$. Then, the spacecraft position and velocity of the planet are matched. The state of the planet at the end of the leg is determined from the arrival time, which is computed from the departure date and the transfer time that the outer loop (Algorithm 2) has previously selected. These constraints are then used to determine optimal values of the inner-loop variables for the leg that minimize the associated impulse.

In the outer loop, the genetic algorithm randomly selects the launch date (in the case of launch events), departure and arrival bodies and the time of flight, retaining only those segments that are feasible. The algorithm then loops over them using the gradient solver to optimize each individual leg, selecting the best combination of segments based on total impulse and total transfer time.

The two primary tasks of the outer-loop heuristic algorithm are to optimize the discrete design variables related to the flyby sequence and to function as an automatic initial guess generator for the inner loop's gradient-based solver. Specifically, the outer loop provides initial guesses of flight times ($\bar{T}$) for each interplanetary leg and launch date ($\bar{t}_0$) based on the lower and upper bounds provided by the user (see Tables 2 and 4). It is worth noting that the number of variables to optimize changes depending on the flyby count for each candidate trajectory. As in MOLTO-IT (Morante et al. 2019), we rely on a method based on null values, known as the null-gene technique, as described in Yam et al. (2011) and Gad and Abdelkhalik (2011). The null values will be skipped by the genetic algorithm when parsing a decision vector, it will only use non-null parameters to build a trajectory. The maximum number of permitted flybys determines the decision vector's maximum dimension.

Each individual in the outer loop population represents a unique inner-loop problem, which entails $N$ legs. Each interplanetary segment is treated as an NLP problem, and involves optimization of a set of parameters defining the leg. The number of NLP decision variables and constraints depends on the terminal event of the current leg. Table 3 lists the inner-loop quantities and their bounds according to the type of event at either end of an interplanetary leg: Launch (L), Flyby (F), or Rendezvous (R). The bounds of the leg transfer time and launch epoch are derived from the outer loop guesses by adding/subtracting a 10% safety margin.

Increasing the margin would improve the chance of convergence for the NLP solver but, according to Morante et al. (2019), could degrade the solution diversity by allowing more individuals to converge to the same time values. The NLP problems are solved with the sequential quadratic programming (SQP) solver implemented in MATLAB's *fmincon* function. The outer loop adds the costs (transfer time and total impulse) of each interplanetary leg (Fig. 6) to compute a fitness estimator for each individual

$$J = [\textstyle\sum_{i=1}^{N} T_i, \sum_{i=1}^{N} \Delta V_i]. \tag{49}$$

For unfeasible trajectories, the inner loop returns a large penalty cost. The genetic algorithm evolves the population in order to minimize the fitness parameter (Eq. 49).

To maintain a fixed variation interval for the angles $\theta_A$ and $\theta_B$, they are redefined as fractions of the total angle traveled during each leg. The initial guess of the in-plane hyperbolic excess velocity angle $\psi_{\infty 0}$ is $\pi/2$ when the second body is an inner planet, $-\pi/2$ for outer planets and 0 for resonant flybys. To improve the convergence of the gradient solver, the orbital elements of the planets are kept constant during each interplanetary segment. At the beginning of the leg, the osculating orbital elements of each body are computed from JPL NAIF-SPICE ephemerides (Acton 1996).

The user only needs to provide the data listed in Table 4 to run Step 1. In Table 4, $t_{0,min}$ and $t_{0,max}$ represent the bounds of the launch epoch, while $n_{\mathrm{fb},min}$ and $n_{\mathrm{fb},max}$ denote the minimum and maximum number of permitted flyby maneuvers. $v_{\infty 0,min}$ and $v_{\infty 0,max}$ stand for the smallest and largest values of the excess hyperbolic speeds, while $h_{\mathrm{fb}min,i}$ is the minimum flyby altitude associated with flyby body $b_i$, and $n$ is the number of candidate flyby bodies. $T_{i,min}$ and $T_{i,max}$ denote the bounds of the transfer time for each leg.

**Table 4.** User inputs for Step 1.

| Description | Values |
|---|---|
| Mission type | Flyby/Rendezvous |
| Launch window | $\{t_{0,\mathrm{min}}, t_{0,\mathrm{max}}\}$ |

| Duration of interplanetary leg | $\{T_{i,\min}, T_{i,\max}\}$ |
|---|---|
| Departure body | $b_0$ |
| Arrival body | $b_f$ |
| Minimum/maximum number of flybys | $\{n_{\text{fb},min}, n_{\text{fb},max}\}$ |
| Flyby bodies | $\{b_1, \ldots, b_n\}$ |
| Minimum flyby radii | $\{h_{\text{fb}min,1}, \ldots, h_{\text{fb}min,n}\}$ |
| Launch $v_{\infty 0}$ | $\{v_{\infty 0,min}, v_{\infty 0,max}\}$ |

The number of optimization variables and constraints in the inner loop increases in comparison with MOLTO-IT (Morante et al. 2019) as a result of incorporating out-of-plane motion. For the launch event, we add a new variable representing the out-of-plane launch excess velocity angle $\eta_0$. When the leg involves a flyby encounter we incorporate the B-Plane angle $\zeta$. In addition, the coefficients $c_2, c_3$ and $c_4$ of the polynomial shaping law that describe the out-of-plane motion must also be determined.

### 4.2    Step 2: High-fidelity optimization

In this step, the candidate trajectories obtained in Step 1 are used as an initial guess for a direct optimization method. A higher-fidelity dynamical model (Eq. 1) is implemented, which also includes changes in mass or additional operational constraints. Thus, the control law is not defined *a priori*; rather, it is determined by a gradient-based solver as a part of the solution, taking propulsion constraints of the LT engine into account. This transforms the multi-objective HOCP into a single-objective MOCP. Adjacent phases in the MOCP are connected by a flyby body, with linkage constraints imposed at the interface to ensure the continuity of mass and position at each boundary. The flyby sequence from Step 1 is used without changes. A direct collocation approach using the Hermite-Simpson collocation scheme (Topputo and Zhang 2014) is used to convert the MOCP into a large-scale parameter optimization problem, which is solved using the NLP solver IPOPT (Biegle and Wächter 2006). In Step 2, only one objective function is permitted, typically the transfer time, the propellant mass, or a weighted combination of both. Additionally, the full set of constraints is enforced.

## 5     3D MOLTO-IT validation: mission to Ceres

This section examines the accuracy and efficiency of the new algorithm when the orbit of the target body is inclined relative to the ecliptic. In order to compare our solution with existing literature, we optimize a rendezvous mission to the asteroid Ceres, whose orbit is inclined 10 degrees. It is not an extreme study case, but results from the literature are available. In particular, this case study compares the performance of the enhanced 3D algorithm against the original MOLTO-IT.

Numerous researchers have computed LT optimal trajectories to Ceres using a variety of strategies, including Sauer (1997), Petropoulos and Longuski (2004), Roa (2017) and Morante et al. (2019). Sauer (1997) found a propellant-efficient route to Ceres via a Mars flyby with a launch date in 2003. His solution often appears in the literature as a benchmark. In 2002, Petropoulos and Longuski (2004) developed the Satellite Tour Design Program-Low Thrust Gravity Assist (STOUR-LGTA) tool based on a two-dimensional exponential sinusoid model. They used it to perform a broad search of optimal rendezvous trajectories to Ceres, all including a flyby with Mars. The approximate solutions from STOUR-LTGA were refined with the direct optimization program GALLOP. Roa (2017) designed an Earth-Mars-Ceres trajectory with a 2003 launch, modeling the trajectory with three-dimensional generalized logarithmic spirals. Recently, Morante et al. (2019) developed the MOLTO-IT tool based on planar generalized logarithmic spirals and used it to perform a multi-objective global search of feasible trajectories to Ceres. They followed with high-fidelity optimization using a direct method, keeping the same 2003 launch window.

For the sake of comparison, we adopt the same transfer configuration as in the previous studies. In particular, we use a 2003 launch window and a launch hyperbolic excess speed of 1.6 km/s. 3D MOLTO-IT performs the Step 1 multi-objective global search of feasible trajectories to Ceres via Mars, Earth and Venus flybys. For Step 2, we assume an initial spacecraft mass of 568 kg, with the same propulsion system as the Deep Space 1 mission (Rayman et al. 2000) (following Morante et al. 2019 and Petropoulos and Longuski 2004). We adopted the propulsion system model of Williams and Coverstone-Carroll (1997):

$$P_a = \frac{P_0}{r^2} \left( \frac{\gamma_0 + \frac{\gamma_1}{r} + \frac{\gamma_2}{r^2}}{1 + \gamma_3 r} \right), \tag{50}$$

$$T = c_{T0} + c_{T1} P_a , \tag{51}$$

$$\dot{m} = c_{m0} + c_{m1} P_a . \tag{52}$$

In the above model, $P_0 = 10 \ kW$ is the power produced by the solar arrays at a distance of 1 astronomical unit (au) from the Sun. The available power $P_a$ (Eq. 50) depends on the distance from the Sun $r$, expressed in au. The values of the coefficients used are $\gamma_0 = 1.1063$, $\gamma_1 = 0.1495$ au, $\gamma_2 = -0.2990$ au$^2$, and $\gamma_3 = -0.0432$ au$^{-1}$. Furthermore, $c_{T0} = -1.9137$ mN, $c_{T1} = 36.2420$ mN/kW, $c_{m0} = 0.47556$ mg/s and $c_{m1} = 0.90209$ mg/(s kW). It is assumed that the thruster can only utilize up to 2.6 kW of electrical power and that it must have at least 0.649 kW to function. The summary of the problem definition is provided in Table 4.

**Table 5.** Earth-Ceres problem description.

|  | **Step 1** | **Step 2** |
|---|---|---|
| Mission type | Rendezvous | |
| Launch window | 01/01/2003 – 31/12/2003 | |
| Departure body | Earth | |
| Arrival body | Ceres | |
| Number of flybys | {0,1,2} | Determined by Step 1 |
| Flyby bodies | Venus, Mars, Earth | Determined by Step 1 |
| Minimum flyby altitude | 200 km | |
| Launch $v_{\infty 0}$ | 1.6 km/s | |
| Specific impulse | 3000 s | Eq. 52 |
| Thrust | - | Eq. 51 |
| Launch mass | - | 568 kg |
| Leg transfer time | 100 – 2000 day | |

Running Step 1 with the parameters shown in Table 4 and the maximum number of generations and population size set both to 100, we obtain the Pareto-front plot in Fig. 5. It contains 3 types of trajectories: Earth-Ceres (E-C), Earth-Mars-Ceres (E-M-C) and Earth-Mars-Mars-Ceres (E-M-M-C). The transfer time ranges between 1.2 and 4.9 years, while the propellant mass fraction lies between 25% and 56%. From these, the user can select the solution that best meets the mission requirements. The solutions with the lowest propellant mass fraction feature a single flyby at Mars, i.e., E-M-C trajectories. This is in agreement with existing literature. Note that, in Step 1, the only assumption about the propulsion system is a constant specific impulse (3000 s) used to determine the propellant mass fraction.

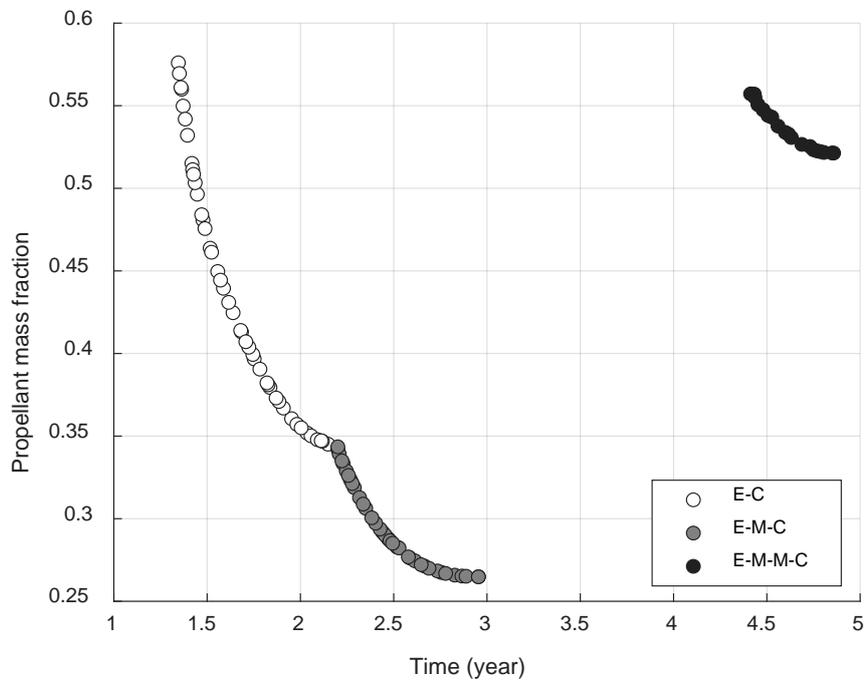

**Fig. 8.** Pareto-Front solution for the Earth-Ceres problem.

We selected the minimum-fuel E-M-C trajectory from Fig. 8 for further optimization in Step 2. Using a tolerance of $10^{-6}$ and 100 nodes to discretize each leg, the NLP solver converged in 106 iterations. Figure 9 shows the difference between the approximate (Step 1) and Step 2 solutions.

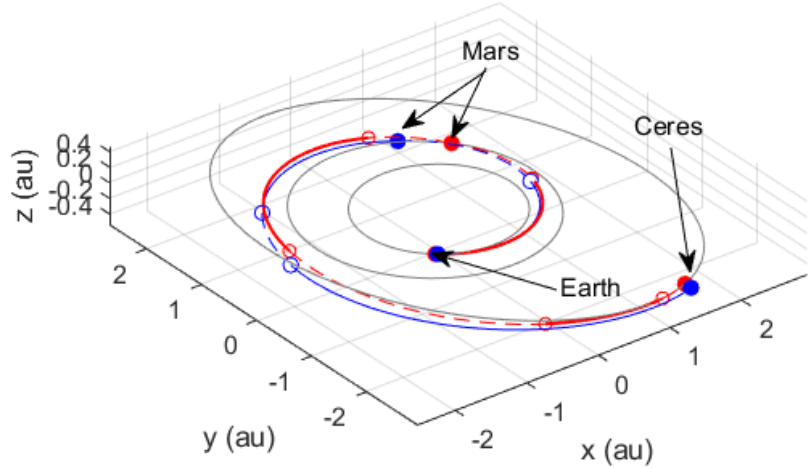

**Fig. 9.** Minimum-propellant E-M-C trajectories obtained after Step 1 (blue) and after Step 2 (red). The continuous and dashed lines represent thrust and coast arcs, respectively.

The comparison against published results is presented in Table 5. The Step 1 and Step 2 solutions in this study are close to those in the literature. Furthermore, the minimum-fuel Step 1 trajectory requires less propellant than those found by Sauer (1997), Roa (2017), and Petropoulos and Longuski (2004), while having a shorter transfer time. This advantage extends to the Step 2 solution when compared against GALLOP (Petropoulos and Longuski 2004).

Comparing the Step 1 solutions of the planar (Morante et al 2019) and 3D versions of MOLTO-IT reveals a higher propellant requirement for the 3D code. However, as explained below, this apparent disadvantage arises because 2D MOLTO-IT disregards completely the out-of-plane motion.

**Table 6.** Comparison of the E-M-C minimum-fuel trajectory with similar published solutions.

|  | 3D MOLTO-IT | | MOLTO-IT Morante 2019 | | Petropoulous & Longuski 2004 | | Other solutions | |
| --- | --- | --- | --- | --- | --- | --- | --- | --- |
|  | **Step 1** | **Step 2** | **Step 1** | **Step 2** | **STOUR-LTGA** | **GALLOP** | **Roa 2017** | **Sauer 1997** |
| Launch date | 13/05 2003 | 12/05 2003 | 2/06 2003 | 13/05 2003 | 6/05 2003 | 6/05 2003 | 21/05 2003 | 8/06 2003 |

| Launch $v_\infty$ (km/s) | 1.60 | 1.60 | 1.60 | 1.60 | 1.60 | 1.60 | 2.00 | 1.60 |
|---|---|---|---|---|---|---|---|---|
| M flyby date | 29/01 2004 | 24/12 2003 | 3/02 2004 | 31/12 2003 | 1/02 2004 | 1/02 2004 | 13/01 2004 | 04/10 2003 |
| Flyby $v_\infty$ (km/s) | 1.96 | 2.34 | 2.40 | 2.16 | 1.43 | 1.92 | 1.27 | - |
| B-Plane angle | 55.1° | 50.2° | 0° | 10.0° | 2.3° | 82.3° | - | - |
| Flyby height (km) | 200 | 200 | 200 | 200 | 5430 | 200 | 2060 | - |
| Arrival date | 07/05 2006 | 06/05 2006 | 18/03 2006 | 26/01 2006 | 12/06 2006 | 9/02 2006 | 16/05 2006 | 06/06 2006 |
| Arrival $v_\infty$ (km/s) | 0.00 | 0.00 | 0.00 | 0.00 | 0.237 | 0.00 | 0.00 | 0.00 |
| Fuel fraction | 0.253 | 0.228 | 0.224 | 0.229 | 0.256 | 0.233 | 0.289 | 0.275 |
| E-M time (day) | 259 | 224 | 216 | 232 | 271 | 271 | 246 | 250 |
| M-C time (day) | 831 | 866 | 774 | 757 | 862 | 739 | 846 | 845 |
| Total time (day) | 1091 | 1091 | 990 | 990 | 1133 | 1110 | 1092 | 1095 |

To extend the comparison between the original and the enhanced tool, we ran the planar code Step 1 with the settings of Table 4. We selected an E-M-C trajectory with the same transfer time as the solution in Fig. 9 as initial guess for Step 2. For both codes, we used the propulsion model in Eqs. 50-52, a grid of 100 nodes per leg and a tolerance of $10^{-6}$ for IPOPT. For this particular problem, Step 2 convergence becomes problematic if the total transfer time is not fixed. Accordingly, we set the total mission duration of the high-fidelity solution equal to the value obtained in Step 1. Note that, while the total length of the mission is constrained, the duration of each leg is free. The results are presented in Table 6.

Both trajectory estimates from Step 1 lead to essentially the same optimal solution after Step 2. However, the 2D guess requires a larger number of iterations to converge (225 vs. 106). It is implied that the 3D MOLTO-IT approximation is closer to the high-fidelity solution. Likewise, the launch and arrival dates from 3D MOLTO-IT differ from the Step 2 solution by one day only. On the other hand, the deviation of the planar code from the optimal solution is 31 days, meaning that the NLP solver in Step 2 corrects the launch date by one month. This is one of the reasons for requiring more iterations to converge.

While the thrust law computed in Step 1 is not used for the subsequent high-fidelity optimizer, it is nonetheless an interesting metric of the quality of the initial guess. Table 6 shows a closer agreement of the thrust impulse values for both legs with 3D MOLTO-IT. For example, the 3D guess differs 7.7% in E-M in-plane impulse from the Step 2 solution, versus 8.7% for the planar approximation. The 3D MOLTO-IT deviation in the M-C out-of-plane impulse is 23%, while there is (obviously) a 100% discrepancy for 2D MOLTO-IT, which does not consider out-of-plane motion. The 3D code relative difference in out-of-plane impulse for the E-M segment seems large (75%), but it is just an effect of the small value (0.44 km/s versus 2.65 km/s for the M-C transfer). Furthermore, the new code approximates the E-M and M-C out-of-plane accelerations to 20.5% and 21.9% accuracy. It also delivers a better estimate of the average in-plane accelerations relative to the Step 2 solution. For the E-M and M-C segments the differences are 3.8% and 12.1%, compared to 3.9% and 15.3% for the planar algorithm. The total impulse estimates by the 3D and 2D codes are 8.62 km/s and 7.31 km/s, versus 8.01 km/s for the accurate solution. Despite its closer agreement in terms of total impulse, the 3D approximation has a higher discrepancy in propellant consumption (10.8% vs. 0.88%). This is caused by the constant impulse assumption made in Step 1. If instead of a variable specific impulse in Step 2 (controlled by Eq. 52) a constant value of 3000 s is prescribed, the propellant mass fraction becomes 0.2464. This translates into a disparity of 10.8% for the 2D guess and 2.8% for its 3D counterpart. On the other hand, the 2D estimate of the Mars $v_\infty$ is better (12.8% disagreement) than the 3D approximation (16.8%). However, the former completely ignores the Mars B Plane angle, while the latter estimates it to an accuracy of 9.7%. Thus, considering the data in Table 6, it is clear that 3D MOLTO-IT outperforms the original code in terms of providing trajectories that are closer to the high-fidelity solution.

**Table 7.** Comparison of results from Step 1 and Step 2 of MOLTO-IT and 3D MOLTO-IT for the minimum-fuel E-M-C trajectory.

|  | MOLTO-IT 2019 | | 3D MOLTO-IT | |
| --- | --- | --- | --- | --- |
|  | Step 1 | Step 2 | Step 1 | Step 2 |
| Launch $v_\infty$ (km/s) | 1.60 | 1.60 | 1.60 | 1.60 |

| Launch mass (kg) | - | 568 | - | 568 |
|---|---|---|---|---|
| Launch date | 12/06 2003 | 12/05 2003 | 13/05 2003 | 12/05 2003 |
| M Flyby date | 11/02 2004 | 22/12 2003 | 29/01 2004 | 24/12 2003 |
| Flyby height (km) | 200 | 200 | 200 | 200 |
| Flyby B-Plane angle | 0º | 50.2º | 55.1º | 50.2º |
| Flyby $v_\infty$ (km/s) | 2.05 | 2.35 | 1.96 | 2.34 |
| E-M in-plane impulse (km/s) | 2.51 | 2.04 | 2.54 | 2.04 |
| M-C in-plane impulse (km/s) | 4.80 | 5.29 | 5.72 | 5.31 |
| E-M off-plane impulse (km/s) | 0.00 | 0.44 | 0.11 | 0.44 |
| M-C off-plane impulse (km/s) | 0.00 | 2.65 | 2.04 | 2.65 |
| E-M time (year) | 0.68 | 0.62 | 0.71 | 0.62 |
| M-C time (year) | 2.31 | 2.37 | 2.28 | 2.37 |
| Total time (year) | 2.99 | 2.99 | 2.99 | 2.99 |
| Mean E-M in-plane accel. ($\mu$m/s$^2$) | 117 | 123 | 115 | 120 |
| Mean E-M off-plane accel. ($\mu$m/s$^2$) | 0.0 | 3.9 | 4.7 | 3.9 |
| Mean M-C in-plane accel. ($\mu$m/s$^2$) | 65.9 | 77.8 | 79.6 | 71.0 |
| Mean M-C off-plane accel. ($\mu$m/s$^2$) | 0.0000 | 35.5 | 27.7 | 35.5 |
| Propellant mass fraction | 0.2209 | 0.229 | 0.253 | 0.229 |
| Arrival date | 02/06 2006 | 06/05 2006 | 07/05 2006 | 06/05 2006 |
| Arrival $v_\infty$ (km/s) | 0.00 | 0.00 | 0.00 | 0.00 |
| NLP iterations | - | 225 | - | 106 |

The good agreement between the 3D MOLTO-IT trajectory and the high-fidelity solution is also demonstrated by the velocity plots in Fig. 10. In all cases, Step 1 and Step 2 solutions follow similar trends with minor deviations. As a crosscheck, the graph also shows the propagated solu-tion using the steering law computed in Step 2. The excellent match between the NLP and propagated solutions demonstrates that the temporal discretization is adequate.

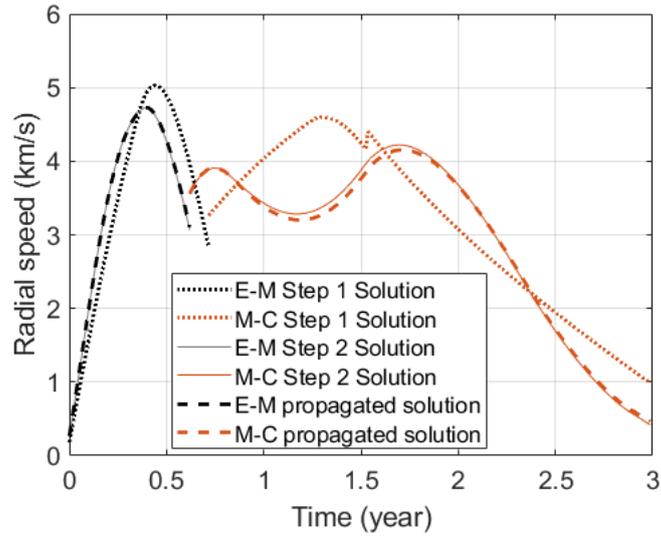

(a) Time history of radial speed.

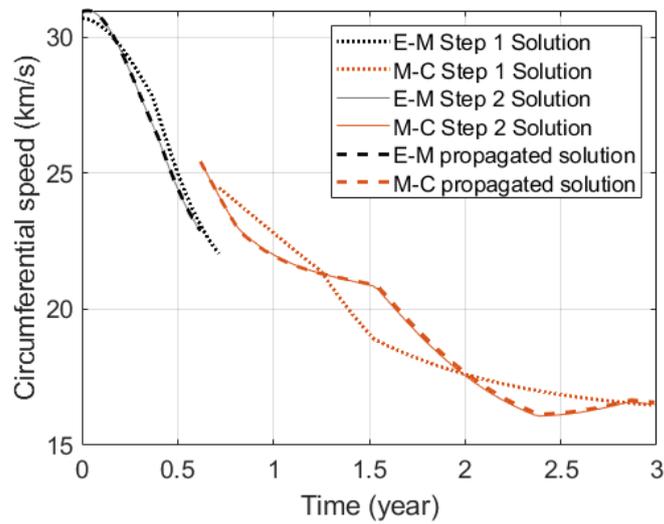

(b) Time history of the circumferential speed.

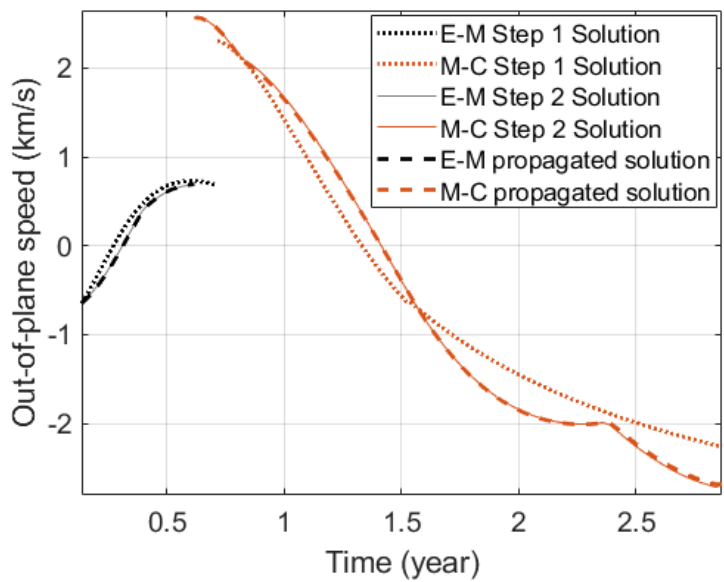

(c) Time history of out-of-plane speed.

**Fig. 10.** Comparison of Step 1 and Step 2 results for the components of velocity.

The agreement between 3D MOLTO-IT Step 1 and Step 2 solutions extends to other trajectories. Figure 11 shows multiple E-M-C solutions with varying transfer time comparing 2D Step 1 solutions with 3D Step 1 solutions and the Step 2 solutions. Here, we compare propellant mass fraction, launch date, Mars flyby date, and arrival date for each trajectory characterized by its time of flight. The Pareto front in Figure 11 (a) allows us to compare the predicted performances of 2D and 3D Step 1, compare to the further optimized Step 2 solutions. It shows a systematic discrepancy between both steps. One of the reasons, as explained above, is the different propulsion system model (constant specific impulse for Step 1 vs. Eqs. 50-52 for Step 2). The 2D Step 1 guesses appear to be closer to the optimal solution than the 3D counterparts due to the former not taking the out-of-plane thrust impulse into account. In other words, the Step 2 is not able to improve the performance of the initial guess provided by the 2D Step 1, because it is initially infeasible, and is in charge of finding the closest feasible solution. From this point of view, 3D Step 1 solutions provide a more conservative picture of the mission performance with a similar trend as the Step 2, and more feasible solutions, in general.

Fig. 11(b-d) shows the launch, flyby and arrival dates. We are using these three variables as a proxy to the computational time invested in Step 2, because it shows how far from the optimal geometry of the Step 2 trajectory is the initial guess., we examine the variation of the launch date with transfer time. There is a close agreement between 3D MOLTO-IT Step 1 and Step 2, with the difference ranging from 0 to 12 days. On the other hand, the discrepancy for the 2D solution is between 29 and 48 days. For the Mars flyby date, Figure 11(c), the difference remains below 36 days for all cases of the 3D guess, and reaches 49 days for the 2D estimate. Additionally, the arrival dates show better agreement of the 3D trajectories with the optimal solutions, with the largest deviation being 12 days. On the other hand, the smallest discrepancy for the 2D guesses is 29 days. The launch, flyby and arrival dates show that 3D MOLTO-IT Step 1 captures the geometry of the optimal solution better, leading to faster convergence in Step 2.

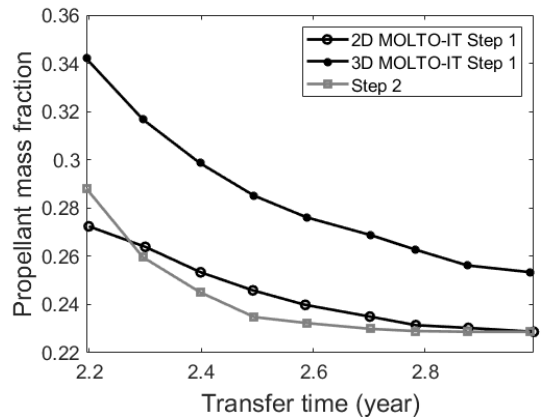

(a)

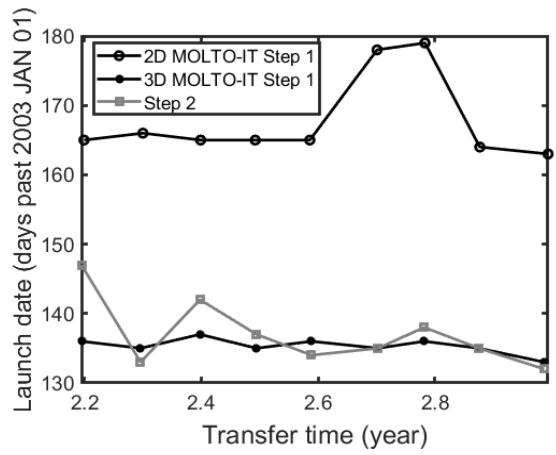

(b)

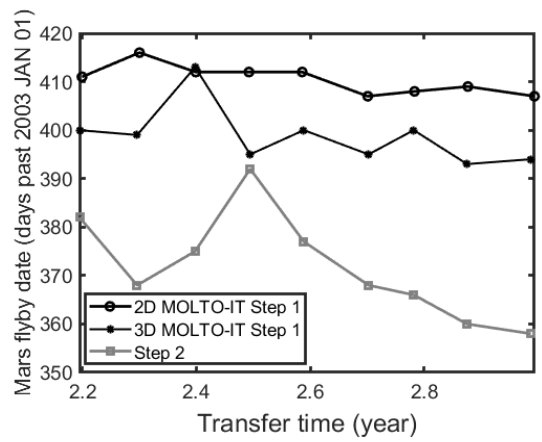

(c)

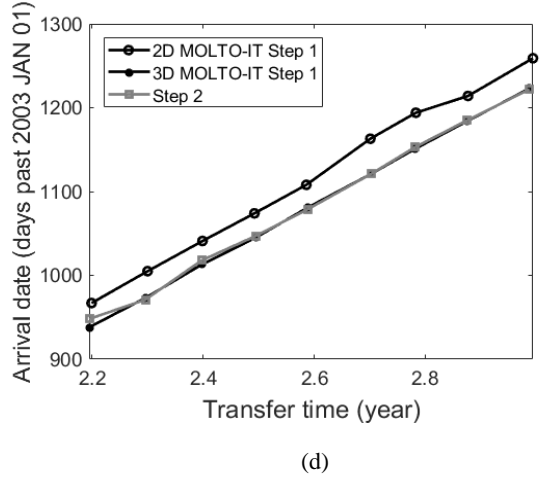

(d)

**Fig. 11.** Comparison between Step 1 and Step 2 solutions for various trajectories.

The improved accuracy in Step 1 comes at the cost of increased computational effort. For example, 3D MOLTO-IT requires 23.4 minutes to complete the Step 1 search for optimal trajectories to Ceres with a maximum of two GA maneuvers with Venus, Mars or Earth using a population size of 100 individuals for 100 generations. For the same settings, 2D MOLTO-IT required 19.4 minutes, i.e., 17% less. This is expected, because 3D MOLTO-IT computes additional variables (around 50% more) and constraints. However, convergence of Step 2 is improved. For the Ceres example, it took 5.6 minutes with the 3D guess vs. 11.9 minutes for the planar solution (a 53% reduction). Thus, considering Steps 1 and 2 together, 3D MOLTO-IT is more efficient than the 2D algorithm.

The run times reported are for a Windows 10 computer (Processor: Intel(R) Xeon(R) Gold 6136 CPU @ 3.00GHz  2.99 GHz, RAM: 64.0 GB) using the MATLAB programming language.

## 6    Conclusions

We presented an automated algorithm for the design of 3D LT multigravity-assist interplanetary trajectories. The work enhances the capability of the MOLTO-IT tool by Morante et al. (2019), which is based on a two-step process not requiring *a priori* knowledge of the flyby sequence. MOLTO-IT assumes a planar surrogate dynamical model for an efficient exploration of the solution space but it is not efficient enough and can potentially fail to provide a feasible solution when the interplanetary trajectory involves visiting high-inclined celestial bodies. In this work, the methodology is therefore expanded to account for out-of-plane motion in the surrogate model.

The first step of the 3D method uses a hybrid heuristic/gradient-based optimizer to find approximate solutions over a large parameter space. It uses predefined trajectory shapes (3D generalized logarithmic spirals for the powered segments) to speed up the computations. We assume two powered segments separated by a coast arc in the case of rendezvous events. For flyby events we use a single thrust segment followed by a ballistic arc. The initial guesses are refined in the second step using a higher fidelity NLP optimizer. The new algorithm improves performance when the mission includes encounters with celestial bodies having moderately or highly inclined orbits, as demonstrated by a rendezvous trajectory to Ceres. The first step of the 3D method model yielded more accurate initial guesses, which improved the convergence of the NLP solver and reduced the total solution time.

# Appendix

### Algorithm 1: Step 1 inner loop

```
Input:
```
$t_0, v_{\infty 0}, \psi_o, \eta_0, \delta, \zeta, T_i, \xi_1, \theta_A, c_{21}, c_{31}, c_{41}, (\xi_2, \theta_B, c_{22}, c_{32}, c_{42})$
```
Output: ΔV
```

```
#Obtain initial state for this leg
if event 1 = launch, then
```
    ```Compute velocity of departure body``` $b_0$`:` $v_{b,0}(t_0)$
    ```Add launcher impulse:``` $v(t_0) = v_{b,0}(t_0) + v_{\infty,0}$
```
else
```
    ```Obtain initial state``` $x(\theta_0)$ ```from leg``` $i-1$
    ```Determine state of flyby planet``` $b_i$`:` $x_{b,i}(t_0)$
    ```Compute the state after flyby```
```
end if
```

```
#Obtain final state
```
    ```Compute final time``` $t_F = t_0 + T_i$
    ```Compute state of planet``` $b_{i+1}$`:` $x_{b,i+1}(t_F)$

    ```if event 2 = flyby```

```
    Propagate thrust arc from θ₀ to θ_A
    Compute the out-of-plane thrust accelera-
tion a_pz
    Propagate coast arc from θ_A to θ_F
  else
    Propagate thrust arc from θ₀ to θ_A
    Compute out-of-plane thrust acceleration
a_pz
    Propagate coast arc from θ_A to θ_B
    Propagate thrust arc from θ_B to θ_F
  end if
```

```
#Enforce constraints
if event 2 = flyby
    cons 1: r‖(θ_F) − r_{b,i+1,‖}(t_F) = 0
    cons 2: r_z(θ_F) − r_{b,i+1,z}(t_F) = 0
    cons 3: t(θ_F) − T_i = 0
    cons 4: a_pz(θ₀) = 0
    cons 5: a_pz((θ₀ + θ_A)/2) = 0
    cons 6: a_pz(θ_A) = 0
else
    cons 7: v‖(θ_F) − v_{b,i+1,‖}(t_F) = 0
    cons 8: ψ(θ_F) − ψ_{b,i+1}(t_F) = 0
    cons 9: v_z(θ_F) − v_{b,i+1,z}(t_F) = 0
end if
```

```
#Compute cost
Obtain ΔV from spirals
```

**Algorithm 2: Step 1 outer loop**

```
Input: t̄₀, b₀, T̄₁, b₁,…, T̄_n, b_n    #Genetic Algorithm
variables
Output: ΔV, ToF

Determine the number of flybys: n_fb
```

```
for i = 0 : n_fb   #Solve the i^th leg
    Departure body = b_i
    Arrival body = b_{i+1}

    #Determine departure and arrival events
    Event 1 = departure   #launch or flyby
    Event 2 = arrival     #flyby or rendezvous

    Use Algorithm 1 to obtain ΔV and T_i for i^th
leg
    if unfeasible then discard
end for

#Compute cost
Compute total ΔV and ToF
```